\documentclass[10pt,a4paper,reqno]{amsart}
\usepackage{amsthm}
\usepackage{amsmath}
\usepackage{amsmath, amssymb}
\usepackage{type1cm}
\usepackage{cite}
\usepackage{bm}
\usepackage{cases}
\usepackage{here}
\usepackage[dvipdfmx]{graphicx}
\theoremstyle{definition}

\newtheorem{prop}{Proposition}
\newtheorem{lemma}{Lemma}
\newtheorem{theorem}{Theorem}
\newtheorem{corollary}{Corollary}
\newcommand{\eqdef}{\stackrel{\mathrm{def}}{=}}
\DeclareMathOperator*{\esssup}{ess\,sup}
\begin{document} 
\title{Inequalities between $L^p$-norms for log-concave distributions}
\author{Tomohiro Nishiyama}
\begin{abstract}
Log-concave distributions include some important distributions such as normal distribution, exponential distribution and so on.
In this note, we show inequalities between two Lp-norms for log-concave distributions on the Euclidean space.
These inequalities are the generalizations of the upper and lower bound of the differential entropy and are also interpreted as a kind of expansion of the inequality between two Lp-norms on the measurable set with finite measure.
\\

\smallskip
\noindent \textbf{Keywords:} log-concave distribution, Lp-norm, moment, maximum entropy, lower bound, differential entropy.
\end{abstract}
\date{}
\maketitle
\bibliographystyle{plain}
\section{Introduction}
A probability density function (pdf) $f:\mathbb{R}^n\rightarrow\mathbb{R}$ is said to be log-concave if
\begin{align}
f(\theta x + (1-\theta) y) \geq f(x)^\theta f(y)^{1-\theta} \nonumber
\end{align}
for each $x,y\in \mathbb{R}^n$ and each $0\leq\theta\leq 1$.
For a random variable $X$ according to a log-concave pdf $f$, the upper and lower bound for the differential entropy is known.
For constants $c_0$ and $c_1$, 
\begin{align} 
\label{entropy_bound}
\log(c_0 \sigma)\leq h(X) \leq \log(c_1\sigma),
\end{align}
where $h(X)\eqdef -\int_{\mathbb{R}} f(x)\log(f(x)) \mathrm{d}x$ denotes the differential entropy for the Lebesgue measure and $\sigma$ denotes the standard deviation of $X$.
$C_1=(2\pi e)^{\frac{1}{2}}$ is a optimal constant\cite{conrad2004probability,zamir1992universal}. Bobkov and Madiman showed $C_0>2^{-\frac{1}{2}}$\cite{bobkov2011entropy} and Marsiglietti and Kostina recently showd tighter bound $C_0=2$ \cite{marsiglietti2018lower}.

Since $h(X)=\lim_{p\rightarrow 1} \frac{1}{1-p}\log \|f\|_p^p$, our motivation is to generalize the inequality (\ref{entropy_bound}) for the $L^p$-norm and the $\alpha$-th moment.

Here, $L^p$-norm\cite{rudin2006real} for the Lebesgue measure is defined as follows.

For $1 \leq p < \infty$,
\begin{align}
\|f\|_p\eqdef {(\int_{\mathbb{R}^n} f(x)^p \mathrm{d}^nx)}^{\frac{1}{p}}.  \nonumber
\end{align}

For $p=\infty$,
\begin{align}
\|f\|_\infty\eqdef \esssup_{x\in \mathbb{R}^n} |f|. \nonumber
\end{align}

In the same way, let us define the $\alpha$-th norm of a random variable $X\sim f$ as follows.
\begin{align}
\sigma_\alpha\eqdef E[|X-E[X]|^\alpha]^{\frac{1}{\alpha}}, \nonumber
\end{align}
where $E[\cdot]$ denotes a expected value for the pdf $f$. From the definition, $\sigma_\alpha^\alpha$ corresponds to the $\alpha$-th moment.

The main purpose in this note is to study relations between $L^p$-norms for log-concave distributions.

For a log-concave pdf $f$, $1\leq p \leq\infty$, $1\leq q \leq\infty$ and $1\leq \alpha < \infty$, 
\begin{align}
\label{ineq_main}
\|f\|_p\leq C_\alpha^{1-\frac{1}{q}}D_\alpha^{1-\frac{1}{p}} \sigma_\alpha^{\frac{1}{p}-\frac{1}{q}}\|f\|_q,
\end{align}
where $C_\alpha, D_\alpha$ are constants which only depend on $\alpha$ (see Theorem 1 in Section 2).
From this inequality we confirm $\|f\|_p \approx \sigma_\alpha^{\frac{1}{p}-1}$ and we derive the upper and lower entropy bound.

Inequality (\ref{ineq_main}) is the similar inequality for the measurable set $\Omega$ with finite measure.
For $1\leq p\leq q\leq \infty$,
\begin{align}
\label{ineq_finite}
\|f\|_{p,\Omega}\leq \mu(\Omega)^{\frac{1}{p}-\frac{1}{q}} \|f\|_{q,\Omega}, 
\end{align}
where $\|f\|_{p,\Omega}\eqdef {(\int_\Omega f(x)^p\mathrm{d}^nx)}^{\frac{1}{p}}$ for $1\leq p <\infty$, $\|f\|_{\infty,\Omega}\eqdef \esssup_{x\in \Omega} |f|$ and $\mu$ is the Lebesgue measure. Since we can interpret $\sigma_\alpha$ and $\mu(\Omega)$ are ``range'' of the regions the pdf $f$ spreads, the inequality (\ref{ineq_main}) is a kind of expansion of (\ref{ineq_finite}).

\section{Main Results}
In the following, the constants $C_\alpha$ and $D_\alpha$ the same for each Theorem, Proposition and Corollary.
\begin{theorem}
\label{th_1dim_norm_inequality}
Let $f$ be a log-concave pdf on $\mathbb{R}$ with finite $\sigma_\alpha$. 

For $1\leq p \leq\infty$, $1\leq q \leq\infty$ and $1\leq \alpha < \infty$,
\begin{align}
\|f\|_p\leq C_\alpha^{1-\frac{1}{q}}D_\alpha^{1-\frac{1}{p}} \sigma_\alpha^{\frac{1}{p}-\frac{1}{q}}\|f\|_q,
\end{align} 
where  $C_\alpha\eqdef\frac{2}{\alpha}\Gamma(\frac{1}{\alpha})(\alpha e)^{\frac{1}{\alpha}}$,  $D_\alpha\eqdef\Gamma(\alpha+1)^{\frac{1}{\alpha}}$ and $\Gamma(x)$ denotes the gamma function.

When $\alpha=2$, the inequality is tighten as 
\begin{align}
\|f\|_p\leq C_2^{1-\frac{1}{q}}\sigma^{\frac{1}{p}-\frac{1}{q}}\|f\|_q, 
\end{align} 
where $\sigma\eqdef \sigma_2$.
\end{theorem}
\begin{corollary}
\label{cor_1dim_norm_inequality}
Let $f$ be a log-concave pdf on $\mathbb{R}$ with finite $\sigma_\alpha$.

For $1\leq p \leq\infty$, $1\leq q \leq\infty$ and $1\leq \alpha < \infty$,
\begin{align}
\label{ineq_both_side}
C_\alpha^{\frac{1}{p}-1}D_\alpha^{\frac{1}{q}-1} \sigma_\alpha^{\frac{1}{p}-\frac{1}{q}}\|f\|_q\leq\|f\|_p\leq C_\alpha^{1-\frac{1}{q}}D_\alpha^{1-\frac{1}{p}} \sigma_\alpha^{\frac{1}{p}-\frac{1}{q}}\|f\|_q.
\end{align} 

When $q=1$, by using $\|f\|_1=1$, (\ref{ineq_both_side}) is simplified as follows.
\begin{align}
\label{ineq_both_side_simple}
{(C_\alpha\sigma_\alpha)}^{\frac{1}{p}-1} \leq \|f\|_p\leq {\biggl(\frac{\sigma_\alpha}{D_\alpha}\biggr)}^{\frac{1}{p}-1}
\end{align}
\end{corollary}
From (\ref{ineq_both_side_simple}), we can confirm $\|f\|_p \approx \sigma_\alpha^{\frac{1}{p}-1}$.
\begin{corollary}
\label{cor_entropy}
Let $X\in\mathbb{R}$ according to a log-concave pdf $f$.

For $1\leq \alpha < \infty$,
\begin{align}
\log\bigl(\frac{\sigma_\alpha}{D_\alpha}\bigr)\leq h(X)\leq \log(C_\alpha\sigma_\alpha).
\end{align} 
\end{corollary}
For $\alpha=2$, this is the same result shown in \cite{bobkov2011entropy}

\begin{prop}
\label{prop_1dim_norm_inequality}
Let $f$ be a symmetric log-concave pdf (that is, $f(x)=f(-x)$) on $\mathbb{R}$ with finite $\sigma_\alpha$.

For $1\leq p \leq\infty$, $1\leq q \leq\infty$ and $1\leq \alpha < \infty$,
\begin{align}
\|f\|_p\leq C_\alpha^{1-\frac{1}{q}} {\biggl(\frac{D_\alpha}{2}\biggr)}^{1-\frac{1}{p}}\sigma_\alpha^{\frac{1}{p}-\frac{1}{q}}\|f\|_q.
\end{align} 
\end{prop}
\begin{theorem}
\label{th_multi_norm_inequality}
Let $f$ be a log-concave pdf on $\mathbb{R}^n$ with finite covariance matrix $\Sigma$.

For $1\leq p \leq\infty$, $1\leq q \leq\infty$ and $n\geq 2$,
\begin{align}
\|f\|_p\leq C(n)^{1-\frac{1}{q}} D(n)^{1-\frac{1}{p}} |\Sigma|^{\frac{1}{2}\bigl(\frac{1}{p}-\frac{1}{q}\bigr)}\|f\|_q,
\end{align} 
where $C(n)\eqdef {(2\pi e)}^{\frac{n}{2}}$, $D(n)\eqdef {\biggl(\frac{n^2e^2}{2\sqrt{2}(n+2)}\biggr)}^{\frac{n}{2}}$ and $|\cdot|$ denotes the determinant of the matrix.
\end{theorem}

\section{Proofs of $L^p$-norm inequalities}
\subsection{Preliminaries for Proofs}
We show some lemmas before the proofs of the main results.
\begin{lemma}
\label{lem_right_1d}
Let $f$ be a pdf on $\mathbb{R}$ with finite $\sigma_\alpha$.

For $1\leq p \leq\infty$ and $0<\alpha < \infty$,
\begin{align}
1\leq C_\alpha^{1-\frac{1}{p}} \sigma_\alpha^{1-\frac{1}{p}} \|f\|_p,
\end{align}
where $C_\alpha\eqdef \frac{2}{\alpha}\Gamma(\frac{1}{\alpha})(\alpha e)^{\frac{1}{\alpha}}$.
\end{lemma}
For $\alpha=2$, tighter bound is shown in \cite{sanchez2011upper}.

\noindent\textbf{Proof.}\\
We prove in the same way as \cite{nishiyama2019p}.
For $\beta > 0$ and a convex function $\phi_t(x)\eqdef \exp\bigl(-\frac{\beta}{t}(x-1)\bigr) $, we consider the following value.
\begin{align}
\label{eq_basic}
V=\int_{\mathbb{R}} f(x)\phi_{p'}(\sigma_{\alpha}^{-\alpha}|x-E[X]|^\alpha) \mathrm{d}x, 
\end{align}
where $p'$ satisfies $\frac{1}{p}+\frac{1}{p'}=1$.

By applying the Jensen's inequality to this equation and using the definition of $\sigma_\alpha $, we get
\begin{align}
\label{eq_jensen}
V\geq \phi_{p'}(\sigma_{\alpha}^{-\alpha}E[|X-E[X]|^\alpha]) =\phi_{p'}(1)=1.
\end{align}

Applying the H\"{o}lder's inequality to (\ref{eq_basic}) yields
\begin{align}
\label{ineq_holder_1dim}
V\leq \|f\|_p{(\int_{\mathbb{R}} \phi_{p'}(\sigma_{\alpha}^{-\alpha}|x-E[X]|^\alpha)^{p'} \mathrm{d}x)}^{\frac{1}{p'}}.
\end{align}
 
By using $\phi_t(x)\eqdef \exp\bigl(-\frac{\beta}{t}(x-1)\bigr) $, we get
\begin{align}
\label{ineq_rhs_1}
\int_{\mathbb{R}} \phi_{p'}(\sigma_{\alpha}^{-\alpha}|x-E[X]|^\alpha)^{p'} \mathrm{d}x=\int_{\mathbb{R}}\exp\bigl(-\beta(\sigma_{\alpha}^{-\alpha}|x-E[X]|^\alpha-1)\bigr)\mathrm{d}x.
\end{align}
Changing from the variable $x$ to $y=\sigma_{\alpha}^{-1}(x-E[X])$ yields
\begin{align}
\label{eq_rhs}
\int_{\mathbb{R}}\exp\bigl(-\beta(\sigma_{\alpha}^{-\alpha}|x-E[X]|^\alpha-1)\bigr)\mathrm{d}x&=2\exp(\beta){\sigma_{\alpha}}\int_0^\infty\exp\bigl(-\beta y^\alpha)\mathrm{d}y \nonumber \\ 
&=\frac{2}{\alpha}\Gamma(\frac{1}{\alpha})\exp(\beta)\beta^{-\frac{1}{\alpha}}\sigma_{\alpha}.
\end{align}
$\exp(\beta)\beta^{-\frac{1}{\alpha}}$ has a minimum value at $\beta=\frac{1}{\alpha}$.
Substituting this condition into (\ref{eq_rhs}) and using (\ref{ineq_holder_1dim}) and (\ref{ineq_rhs_1}) give
\begin{align}
\label{ineq_holder_1dim_2}
V\leq \|f\|_p{(C_\alpha{\sigma_{\alpha}})}^{\frac{1}{p'}},
\end{align}
where $C_\alpha=\frac{2}{\alpha}\Gamma(\frac{1}{\alpha})(\alpha e)^{\frac{1}{\alpha}}$.

By combining (\ref{eq_jensen}),  (\ref{ineq_holder_1dim_2}) and $\frac{1}{p}+\frac{1}{p'}=1$, 
the result follows.
\begin{lemma}
\label{lem_right_multi}
Let $f$ be a pdf on $\mathbb{R}^n$ with finite covariance matrix $\Sigma$.

For $1\leq p \leq\infty$ ,
\begin{align}
1\leq {(C(n){|\Sigma|}^{\frac{1}{2}})}^{1-\frac{1}{p}}\|f\|_p,
\end{align}
where $C(n)\eqdef {(2\pi e)}^{\frac{n}{2}}$.
\end{lemma}

\noindent\textbf{Proof.}\\
We can prove this lemma in the same way as the Lemma \ref{lem_right_1d}.
For a convex function $\phi_t(x)=\exp\bigl(-\frac{\beta}{t}(x-n)\bigr) $, we consider the following value.
\begin{align}
\label{eq_multi_basic}
V=\int_{\mathbb{R}^n} \phi_{p'}((x-E[X])^T\Sigma^{-1}(x-E[X])) \mathrm{d}^nx,
\end{align}
where $p'$ satisfies $\frac{1}{p}+\frac{1}{p'}=1$ and $T$ denotes transpose of a matrix.
By applying the Jensen's inequality to this equation in the same way as Lemma \ref{lem_right_1d}, we get
\begin{align}
\label{eq_multi_jensen}
V\geq \phi_{p'}(n)=1.
\end{align}
Next, by applying the H\"{o}lder's inequality to (\ref{eq_multi_basic}), we have
\begin{align}
\label{eq_multi_holder_1}
V\leq \|f\|_p{\biggl(\int_{\mathbb{R}^n} \phi_{p'}\bigl((x-E[X])^T\Sigma^{-1}(x-E[X])\bigr)^{p'} \mathrm{d}^nx\biggr)}^{\frac{1}{p'}}.
\end{align}
Substituting $\phi_{p'}(x)=\exp\bigl(-\frac{\beta}{p'}(x-n)\bigr)$ into this inequality gives
\begin{align}
\label{eq_multi_integral_1}
\int_{\mathbb{R}^n} \phi_{p'}\bigl((x-E[X])^T\Sigma^{-1}(x-E[X])\bigr)^{p'} \mathrm{d}^nx=\exp(\beta n)\int_{\mathbb{R}^n}\exp\bigl(-\beta(x-E[X])^T\Sigma^{-1}(x-E[X])\bigr)\mathrm{d}^nx.
\end{align}
Changing the variable from $x$ to $y=\Sigma^{-\frac{1}{2}}(x-E[X])$ gives
\begin{align}
\label{eq_multi_rhs}
\exp(\beta n)\int_{\mathbb{R}^n}\exp\bigl(-\beta(x-E[X])^T\Sigma^{-1}(x-E[X])\bigr)\mathrm{d}^nx=\exp(\beta n){\biggl(\frac{\pi}{\beta}\biggr)}^{\frac{n}{2}}{|\Sigma|}^{\frac{1}{2}}.
\end{align}
$\exp(\beta)\beta^{-\frac{1}{2}}$ has a minimum value at $\beta=\frac{1}{2}$.
Substituting this condition into (\ref{eq_multi_rhs}) and combining (\ref{eq_multi_holder_1}) and (\ref{eq_multi_integral_1}) give
\begin{align}
\label{eq_multi_holder_2}
V\leq \|f\|_p{(C(n){|\Sigma|}^{\frac{1}{2}})}^{\frac{1}{p'}},
\end{align}
where $C(n)={(2\pi e )}^{\frac{n}{2}}$.
By combining (\ref{eq_multi_jensen}),  (\ref{eq_multi_holder_2}) and $\frac{1}{p}+\frac{1}{p'}=1$, the result follows.

\begin{lemma}
\label{lem_left}
Let $f$ be a pdf on $\mathbb{R}^n$ with $\|f\|_\infty < \infty$.

For $1\leq p \leq\infty$,
\begin{align}
\|f\|_p \|f\|_\infty^{\frac{1}{p}-1}\leq 1.
\end{align}
\end{lemma}
\noindent\textbf{Proof.}\\
If $p=\infty$, equality holds.

When $1\leq p < \infty$,
\begin{align}
\|f\|_p^p=\int_{\mathbb{R}^n} f(x)^p\mathrm{d}^nx=\int_{\mathbb{R}^n} f(x)^{p-1} f(x)\mathrm{d}^nx
\leq \|f\|_\infty^{p-1}.
\end{align}
Hence, the result follows.
\begin{lemma}
\label{lem_square_infty}
Let $f$ be a log-concave pdf on $\mathbb{R}^n$.

Then,
\begin{align}
\|f\|_\infty \leq 2^n\|f\|_2^2.
\end{align}
\end{lemma}
\noindent\textbf{Proof.}\\
By the definition of log-concavity, for any $x,y\in\mathbb{R}^n$,
\begin{align}
f\bigl(\frac{x}{2}+\frac{y}{2}\bigr)^2\geq f(x)f(y).
\end{align}
Integrating with respect to $x$ and using $\|f\|_1=1$, we get
\begin{align}
2^n\|f\|_2^2\geq f(y).
\end{align}
Optimizing over $y$ yields
\begin{align}
\label{ineq_square}
2^n\|f\|_2^2\geq \|f\|_\infty.
\end{align}

\begin{lemma}
\label{lem_ineq_max_1d}
Let $f$ be a log-concave pdf on $\mathbb{R}$ with finite $\sigma_\alpha$.

For $1\leq \alpha < \infty$,
\begin{align}
\|f\|_\infty \leq \frac{\Gamma(\alpha+1)^{\frac{1}{\alpha}}}{\sigma_\alpha}.
\end{align}

When $\alpha=2$, the inequality is tighten as 
\begin{align}
\label{ineq_tighter}
\|f\|_\infty \leq \frac{1}{\sigma}.
\end{align}

\end{lemma}
\noindent\textbf{Proof.}\\
For a symmetric log-concave random variable $Z\in\mathbb{R}\sim F$, it was shown that $F(Z)$ satisfies the following inequality (see \cite{marsiglietti2018lower}).
\begin{align}
\label{ineq_sym_1d}
F(0)\leq \frac{\Gamma(\alpha+1)^{\frac{1}{\alpha}}}{2E[|Z|^\alpha]^{\frac{1}{\alpha}}}.
\end{align}
Next, we consider random variables $X$ and $Y$ according to $f$.
When $X$ and $Y$ are independent, $X-Y\sim f_{X-Y}$ is symmetric and log-concave. Furthermore, $f_{X-Y}$ satisfies
\begin{align}
\label{eq_x_y}
f_{X-Y}(0)=\|f\|_2^2.
\end{align}
By combining this equation and (\ref{ineq_sym_1d}), we get
\begin{align}
\label{ineq_square_moment}
\|f\|_2^2\leq \frac{\Gamma(\alpha+1)^{\frac{1}{\alpha}}}{2E[|X-Y|^\alpha]^{\frac{1}{\alpha}}}.
\end{align}
From the Jensen's inequality, we get $E[|X-Y|^\alpha]\geq E[|X-E[X]|^\alpha]=\sigma_\alpha^\alpha$.
Hence, we get
\begin{align}
\|f\|_2^2\leq \frac{\Gamma(\alpha+1)^{\frac{1}{\alpha}}}{2\sigma_\alpha}.
\end{align}
Combining this inequality and Lemma \ref{lem_square_infty} for $n=1$ yields the desired result.
When $\alpha=2$, we get $E[|X-Y|^2]^{\frac{1}{2}}=\sqrt{2}\sigma$.
By combining this equation and (\ref{ineq_square_moment}), we can prove (\ref{ineq_tighter}) in the same way.
 
\begin{lemma}
\label{lem_ineq_max_multi}
Let $f$ be a log-concave pdf on $\mathbb{R}^n$ with a finite covariance matrix $\Sigma$.

For $n\geq 2$,
\begin{align}
\label{ineq_sym_multi}
\|f\|_\infty \leq \frac{D(n)}{|\Sigma|^{\frac{1}{2}}}, 
\end{align}
where $D(n)={\biggl(\frac{n^2e^2}{2\sqrt{2}(n+2)}\biggr)}^{\frac{n}{2}}$
\end{lemma}
\noindent\textbf{Proof.}\\
For a symmetric log-concave random vector $Z\in\mathbb{R}^n\sim F$ with covariance matrix $\Sigma_Z$, it was shown that $\tilde{Z}\eqdef \Sigma_Z^{-\frac{1}{2}}Z\sim \tilde{F}$ satisfies the following inequality (in detail, see the proof of Theorem 4 in \cite{marsiglietti2018lower}).
\begin{align}
\tilde{F}(0)\leq 2^{-\frac{n}{2}}D(n).
\end{align}
Since $\tilde{F}(\tilde{Z})=F(Z)|\Sigma_Z|^{\frac{1}{2}}$, we get 
\begin{align}
\label{ineq_sym_multi}
F(0)\leq \frac{2^{-\frac{n}{2}}D(n)}{|\Sigma_Z|^{\frac{1}{2}}}.
\end{align}

Next, we consider random vectors $X$ and $Y$ according to $f$ as well as Lemma \ref{lem_ineq_max_1d}.
When $X$ and $Y$ are independent, $X-Y\sim f_{X-Y}$ is symmetric and log-concave and the covariance matrix of $X-Y$ satisfies $\Sigma_{X-Y}=2\Sigma$ and $f_{X-Y}(0)=\|f\|_2^2$.

By combining these equations and (\ref{ineq_sym_multi}), we get
\begin{align}
\|f\|_2^2\leq  \frac{D(n)}{2^n|\Sigma|^{\frac{1}{2}}}.
\end{align}

Combining this inequality and Lemma \ref{lem_square_infty} yields the desired result.
\subsection{Proofs of Main Results}
We prove inequalities of main results by applying Lemmas shown in the previous subsection.

\noindent\textbf{Proof of Theorem \ref{th_1dim_norm_inequality}.}\\
Combining Lemma \ref{lem_right_1d} for $q$ and Lemma \ref{lem_left} for $p$ yields
\begin{align}
\|f\|_p \|f\|_\infty^{\frac{1}{p}-1}\leq 1\leq C_\alpha^{1-\frac{1}{q}}\sigma_\alpha^{1-\frac{1}{q}} \|f\|_q.
\end{align}
Hence, we get 
\begin{align}
\label{eq_right_1d}
\|f\|_p \leq \|f\|_\infty^{1-\frac{1}{p}}C_\alpha^{1-\frac{1}{q}}\sigma_\alpha^{1-\frac{1}{q}} \|f\|_q.
\end{align}
From Lemma \ref{lem_ineq_max_1d}, the result follows.
We can prove the tighter inequality for $\alpha=2$ in the same way.

\noindent\textbf{Proof of Corollary \ref{cor_1dim_norm_inequality}.}\\
Exchanging $p$ and $q$ in Theorem \ref{th_1dim_norm_inequality} yields
\begin{align}
C_\alpha^{\frac{1}{p}-1}D_\alpha^{\frac{1}{q}-1} \sigma_\alpha^{\frac{1}{p}-\frac{1}{q}}\|f\|_q\leq\|f\|_p.
\end{align}
By combining this inequality and Theorem \ref{th_1dim_norm_inequality}, the result follows.

\noindent\textbf{Proof of Corollary \ref{cor_entropy}.}\\
From (\ref{ineq_both_side_simple}), we have
\begin{align}
\log\bigl(\frac{\sigma_\alpha}{D_\alpha}\bigr)\leq h_p(X)\leq \log(C_\alpha\sigma_\alpha).
\end{align}
where $h_p(X)\eqdef \frac{1}{1-p}\log\|f\|_p^p$ is the R\'{e}nyi entropy \cite{renyi1961measures}.
In the limit $p\downarrow 1$, the result follows.

\noindent\textbf{Proof of Proposition \ref{prop_1dim_norm_inequality}.}\\
By combining Lemma \ref{lem_right_1d}, Lemma \ref{lem_left} and (\ref{ineq_sym_1d}), we can prove in the same way as Theorem \ref{th_1dim_norm_inequality}.

\noindent\textbf{Proof of Theorem \ref{th_multi_norm_inequality}.}\\
By combining Lemma \ref{lem_right_multi}, Lemma \ref{lem_left} and Lemma \ref{lem_ineq_max_multi}, we can prove in the same way as Theorem \ref{th_1dim_norm_inequality}.

 
\section{Conclusion}
For the log-concave pdf $f$ and the $\alpha$-th norm of the random variable $X\sim f$, we have confirmed that $\|f\|_p\approx\sigma_\alpha^{\frac{1}{p}-1}$ for $1\leq p\leq\infty$ and have shown inequalities between two $L^p$-norms. We have also shown these inequalities are the generalizations of the upper and lower bound of the differential entropy.

It is the future work to confirm whether similar inequalities hold or not for $p,q\leq 1$.

\bibliography{reference_v1}

\begin{thebibliography}{1}

\bibitem{bobkov2011entropy}
Sergey Bobkov and Mokshay Madiman.
\newblock The entropy per coordinate of a random vector is highly constrained
  under convexity conditions.
\newblock {\em IEEE Transactions on Information Theory}, 57(8):4940--4954,
  2011.

\bibitem{conrad2004probability}
Keith Conrad.
\newblock Probability distributions and maximum entropy.
\newblock {\em Entropy}, 6(452):10, 2004.

\bibitem{marsiglietti2018lower}
Arnaud Marsiglietti and Victoria Kostina.
\newblock A lower bound on the differential entropy of log-concave random
  vectors with applications.
\newblock {\em Entropy}, 20(3):185, 2018.

\bibitem{nishiyama2019p}
Tomohiro Nishiyama.
\newblock $l^p$-norm inequality using q-moment and its applications.
\newblock {\em arXiv preprint arXiv:1902.01021}, 2019.

\bibitem{renyi1961measures}
Alfr{\'e}d R{\'e}nyi et~al.
\newblock On measures of entropy and information.
\newblock In {\em Proceedings of the Fourth Berkeley Symposium on Mathematical
  Statistics and Probability, Volume 1: Contributions to the Theory of
  Statistics}. The Regents of the University of California, 1961.

\bibitem{rudin2006real}
Walter Rudin.
\newblock {\em Real and complex analysis}.
\newblock Tata McGraw-Hill Education, 2006.

\bibitem{sanchez2011upper}
Pablo Sanchez-Moreno, Steeve Zozor, and Jesus~S Dehesa.
\newblock Upper bounds on shannon and r{\'e}nyi entropies for central
  potentials.
\newblock {\em Journal of Mathematical Physics}, 52(2):022105, 2011.

\bibitem{zamir1992universal}
Ram Zamir and Meir Feder.
\newblock On universal quantization by randomized uniform/lattice quantizers.
\newblock {\em IEEE Transactions on Information Theory}, 38(2):428--436, 1992.

\end{thebibliography}
\end{document}